\documentclass[preprint,12pt]{elsarticle}

\usepackage{amssymb}
\usepackage{amsmath}
\usepackage{amsthm}

\newtheorem{thm}{Theorem}[section]

\usepackage{graphicx}
\usepackage{xcolor}
\usepackage{tikz}
\graphicspath{ {./images/}}
\usepackage{caption}
\usepackage{subcaption}
\usepackage{float}
\usepackage{hyperref}
\usepackage[linesnumbered,ruled,vlined]{algorithm2e}

\usepackage{algorithmic}
\usepackage{enumitem}

\usepackage{multirow}
\usepackage{natbib}

\journal{Applied Numerical Mathematics}

\begin{document}

\begin{frontmatter}



\title{Exponential Convergence of Deep Composite Polynomial Approximation for Cusp-Type Functions}


\author[uchic]{Kingsley Yeon}

\affiliation[uchic]{organization={Department of Statistics and CCAM, University of Chicago},
            city={Chicago},
            postcode={60637}, 
            state={IL},
            country={USA}}

\author[zbmath]{Steven B. Damelin}

\affiliation[zbmath]{organization={Department of Mathematics, ZBMATH-OPEN},
            addressline={Hermann-von-Helmholtz-Platz 1}, 
            city={Eggenstein-Leopoldshafen},
            postcode={76344}, 
            country={Germany}}

\author[hebrew]{Michael Werman}

\affiliation[hebrew]{organization={Department of Computer Science, The Hebrew University},
            city={Jerusalem},
            postcode={91904}, 
            country={Israel}}

\begin{abstract}
We investigate  deep composite polynomial approximations of continuous but non-differentiable
functions with algebraic cusp singularities.
The functions in focus  consist of 
finitely many cusp terms of the form $|x-a_j|^{\alpha_j}$ with rational exponents
$\alpha_j\in(0,1)$ on a real-analytic background.
We propose a constructive approximation scheme that combines a division-free polynomial
iteration for fractional powers with an outer layer for the analytic polynomial fitting.
Our main result shows that this composite structure achieves exponential convergence in
the  the  number of scalar coefficients in the
inner and outer polynomial layers.
Specifically, the $L^p([-1,1])$ approximation error,  decays
exponentially with respect to the parameter budget, in contrast to the algebraic rates
obtained by classical single-layer polynomial approximation for cusp-type functions.
Numerical experiments for both single and multiple cusp configurations confirm the
theoretical rates and demonstrate the parameter efficiency of deep composite polynomial
constructions.
\end{abstract}







\end{frontmatter}



\section{Introduction}

A central problem of approximation theory is finding for a given a function $f$ from a normed function space $X$, a function $g$ in $X$, norm close to $f$. The distance of $f$ to $g$ is then given by the norm $||f-g||_{X}$ and one has to then, measure if possible, how good this approximation is.
This motivates the following. 
\medskip

One form of the  classical theorem of Weierstrass asserts that every continuous function $f:[-1,1]\to \mathbb R$ can be arbitrarily and uniformly approximated by polynomials. More precisely, the theorem says the following. Given any small tolerance $\varepsilon>0$ and continuous $f:[-1,1]\to \mathbb R$, there is a polynomial $R$ with $|f(x)-R(x)|<\varepsilon$. Here, $-1\leq x\leq 1$. That is, the polynomials are dense in the space of real valued  continuous functions on $[-1,1]$. The theorem, however, does not say how high the degree of a polynomial should be so that it can approximate a given function with a preassigned degree of accuracy, and it seems natural to guess that the smoother the function, the less would be the degree the polynomial to accomplish the job. In this paper we study convergence of deep algebraic polynomial approximation to continuous but non-differentiable real valued functions on $[-1,1]$, see Theorem~\ref{thm:deep} below. \footnote{$C([-1,1])$ will henceforth denote the space of 
real valued continuous functions $f:[-1,1]\to \mathbb R$.}
To set the scene for our work, we need to recall some classic approximation results for various real valued Banach function classes $X$ where the functions have domain $[-1,1]$ and are approximated by algebraic polynomials of given degree. The approximation will be done via an interplay of degree of the polynomial space and level of smoothness of the approximation function space $X$. In this regard, let $\Pi_n$ denote the class of all algebraic polynomials of degree at most $n\geq 1$. Then for $f\in X$, define $E_{n}(f)_{X}={\rm min}_{P\in \Pi_n}||f-P||_{X}$, to be the error in best approximation of a fixed function $f$ by elements of $\Pi_n$ when it exists. \footnote{In this paper, we choose to exclude in our study approximation by trigonometric polynomials $\frac{a_0}{2}+\sum_{k=1}^n
(a_k\cos(kx)+b_k\sin(kx))$ (for real numbers $a_i$ and $b_i$) of degree at most $n$ on the unit circle in our analysis, the latter understood as $\mathbb R$ with the identification of points modulo $2\pi$. We defer this to a future paper where we will translate naturally from $[-1,1]$ by appropriate substitutions.} 

The purpose of the present work is to demonstrate that deep composite polynomial constructions
can substantially outperform classical single-layer polynomial approximation for a broad class
of continuous but non-differentiable functions with algebraic cusp singularities.
Our approach combines a division-free polynomial iteration that resolves fractional power
singularities with an outer analytic polynomial fit, yielding a constructive approximation scheme
with provable exponential convergence in the parameter count.
The classical bounds recalled therefore serve as a natural baseline against which the
advantages of deep composite polynomial approximation can be assessed.

\section{Degree of Approximation}
\subsection{Some Notation}
We need some context for measurements both for the degree of approximation and for the smoothness of the function in question. We will work with the following Banach spaces $X$ of real valued functions with domain $[-1,1]$ with the given norms: $L_p$ norms given by $||f||_{p}:=||f||_{p}
:=\left(\int_{-1}^{1}|f(x)|^{p}dx\right)^{1/p}$, 
for $1\leq p<\infty$ and the $L_{\infty}$ norm given by $||f||_{\infty}:=||f||_{\infty}
={\rm max}_{[-1,1]}|f(x)|$ both when finite. For $0<p<1$, we interpret $L_p$ as a \emph{quasi-norm}. Given $X$ as above, we will also work with the following Sobolev function space defined as follows. Let $1\leq p\leq \infty$, $r=1,2,..$ and denote by $W_p^{r}(X)$, the space of real valued functions $f:[-1,1]\to \mathbb R$ for which $f^{(r-1)}$ is absolutely continuous and $f^{(r)}$ is in $X$. The space $W_p^{r}(X)$ is endowed with the norm: $||f||_{W_p^{r}(X)}:=||f||_{X}+||f^{(r)}||_{X}.$
The choice $X=L_{p}$ for $1\leq p<\infty$ yields the Sobolev spaces $W_p^{r}(X):=W_p^{r}$. The choice $X=C[-1.1]$, gives the space $C^{(r)}$ of $r$-times continuously differentiable functions on $[-1,1]$ and the space $C^{\infty}$ consists of infinitely differentiable real valued functions
on $[-1,1]$. For the function spaces above, we will write $E_n(f)_{X}:=E_n(f)_p$.
\footnote{It is known that $E_n(f)_p$ exists and is well defined.} 
\medskip
Throughout, $C,c,C_1,...$ will denote positive constants depending on different parameters  which we will indicate as needed. The same symbol may denote a different constant at any time. The context will be clear.
Letters such as $f,g,P,H..$ will denote functions or classes of functions of various kinds. The same symbol may denote a different function at any time. The context will be clear.

\subsection{Moduli of continuity}

It is well known that obtaining upper bounds for $E_n(f)_{p}$ using levels of differentiability turns out to be too crude. To this end, we recall needed moduli of continuity. 
\medskip

We recall the forward difference operators of a function $f:[-1,1]\to \mathbb R$. Let $h\in \mathbb R$,
$r\geq 1$ be an integer. Set $\Delta_{h}^{1}f(x):=\Delta_{h}f(x):=f(x+h)-f(x)$ and inductively define
$\Delta_{h}^{r}:=\Delta_{h}(\Delta_{h}^{r-1}f(x)),\, r\geq 2$. Set $A_{rh}=[-1,1-rh]$.
Then for $f\in L_p$, $1\leq p\leq \infty$, set for $t\geq 0$, 
$w_r(f,t)_p:={\rm sup}_{0<h\leq t}||\Delta_{h}^{r}f(.)||_p(A_{rh})$ and when $r=1$, write 
$w_r(f,t)=w(f,t)$.
\medskip

\subsection{ Rates of approximation by algebraic polynomials in degree and via smoothness: Upper bounds}

There are many classical results on upper bounds for $E_n(f)_{X}$ in degree and via smoothness. We state the following two from
\cite{devorelorentz} and \cite{mhaskar}.

\medskip

\begin{thm}
For functions $f\in W_p^{r}$, $r=0,1,2,...$, $1\leq p\leq \infty$, the error of algebraic polynomial approximation in degree and smoothness satisfies:
\[
E_n(f)_p\leq Cn^{-r}w(f^{(r)},1/n)_p,\, n>r.
\]
Here $C$ is independent of $r,n$.
In particular, for $f\in W_p^1$,
\[
E_n(f)_p\leq \frac{\pi}{2n+2}||f'||_p,
\]
\end{thm}

\begin{thm}
For each $r=1,2,...$, there is a constant $C$ depending on $r$ (not on $n$) such that for 
$f\in L_p$, $1\leq p\leq \infty$, the error of algebraic polynomial approximation in degree and smoothness satisfies
\[
E_n(f)_p\leq Cw_r(f,1/n)_p,\, n\geq r.
\]
In particular,
\[
E_n(f)_p\leq Cn^{-r}||f^{(r)}||_p
\]
with $C$ independent of $n,r$.
\end{thm}

\section{Motivation}

Our work in this paper, builds on a line of research on deep composite polynomial approximation.
The original motivation appeared in \cite{yeon2025deep}, which showed that composite polynomials
formed by nesting low-degree polynomials can approximate non-smooth functions such as the absolute
value with exponential convergence in the number of free parameters.
This framework was extended to weighted approximation on unbounded domains in
\cite{yeon2025uniform}, with one-sided weights providing uniform control of functions with distinct two-sided behavior, i.e., decay to zero on one side and growth to infinity on the other.
In practice, composite polynomials have been shown to be optimally efficient in terms of the number of non-scalar multiplications, and have been used for the efficient evaluation of matrix functions \cite{jarlebring2025polynomial,yeon2025bolt}.
\medskip

In this paper we study convergence of deep algebraic polynomial approximation to continuous but non-differentiable functions. Algorithms that utilize only addition and multiplication are crucial for homomorphic encryption, a powerful cryptographic technique that enables computations on encrypted data without requiring access to the plaintext \cite{rivest1978ondata,EMS_homomorphic}. This characteristic is particularly significant in scenarios where privacy and data security are paramount, such as in cloud computing and secure data sharing \cite{researchgate_homomorphic,Acar2018ASO}.
In  applied fields, such as physics and engineering, cusps are significant for understanding dynamic systems and phase transitions. For example, in fluid dynamics, cusps can represent critical points where the behavior of fluids changes dramatically, impacting the design of systems such as pipelines or aircraft. Similarly, in materials science, the study of cusps in stress-strain curves helps in identifying failure points in materials under different loads, which is critical for ensuring safety and reliability in engineering applications, \cite{Callister_Rethwisch_2018}.

In computer vision, detecting cusps in image contours is essential for object recognition and shape analysis, allowing for more accurate interpretations of visual data, \cite{Szeliski_2010}.  Locating  the presence and nature of cusps also allows for the development of algorithms to  smooth curves or surfaces without losing essential features, thus enhancing visual realism in digital models, \cite{koenderink1990solid}.

\section{Deep composite polynomial approximation for cusp--type continuous functions}

The functions we approximate have a finite number of algebraic cusp singularities composed with analytic envelopes. Concretely, fix points $a_1,\dots,a_M\in[-1,1]$ and rational exponents
\[
\alpha_j=\frac{r_j}{s_j}\in(0,1),\qquad r_j,s_j\in\mathbb{N},\ s_j\ge 2.
\]
Let $H:[-1,1]\to\mathbb{R}$ be real analytic in a complex neighborhood of $[-1,1]$ and, for each $j$, let $h_j:[0,1]\to\mathbb{R}$ be real analytic in a complex neighborhood of $[0,1]$. Define the class
\[
\mathcal{F}
=\biggl\{\,f(x)=H(x)+\sum_{j=1}^M h_j\!\bigl(|x-a_j|^{\alpha_j}\bigr)\;:\; x\in[-1,1]\,\biggr\}.
\]
Each $f\in\mathcal{F}$ is continuous on $[-1,1]$ and may fail to be differentiable at the points $a_j$ with algebraic cusp profile $|x-a_j|^{\alpha_j}$.


\subsection{Division--free inner iteration and outer analytic fit}

Fix a single cusp location $a\in[-1,1]$ and a single rational exponent $\alpha=r/s\in(0,1)$ with $r,s\in\mathbb{N}$, $s\ge 2$. For $t\in[0,1]$ consider the polynomially defined coupled iteration for the $s$th root \cite{blanchard2023newton}:
\begin{equation}\label{eq:divfree}
\boxed{\quad
\begin{aligned}
y_{k+1}(t) &= y_k(t)\Bigl(2- s\,g_k(t)^{\,s-1}\,y_k(t)\Bigr),\\[2mm]
g_{k+1}(t) &= g_k(t) - y_{k+1}(t)\,\Bigl(g_k(t)^{\,s}- t\Bigr),
\end{aligned}
\qquad g_0\equiv 1,\ \ y_0\equiv \frac1s. \quad}
\end{equation}
Every update uses only additions and multiplications, hence by induction $g_k,y_k$ are polynomials in $t$. Set
\[
\phi_k(t):=g_k(t)^{\,r}\qquad\text{(a polynomial in $t$)}.
\]
The deep approximant we will use for one cusp term is
\[
G_{m,k}(x):=H(x)+P_m\!\left(\phi_k(|x-a|)\right),
\]
where $P_m$ is a polynomial of degree $m$ that approximates the analytic envelope $h$ on $[0,1]$. The effective parameter count of $G_{m,k}$ is
\[
N=\underbrace{(m+1)}_{\text{outer coefficients}}+\underbrace{O(k)}_{\text{inner layers}}
\asymp m+k,
\]
since each inner update \eqref{eq:divfree} introduces $O(1)$ new scalar coefficients and the outer polynomial contributes $m+1$ coefficients.

\begin{thm}[Exponential--in--parameters deep approximation]\label{thm:deep}
Let $0<p<\infty$. For the class $\mathcal{F}$ defined above, there exist constants $C,c>0$ depending only on $p$ with the following property. For every $f\in\mathcal{F}$ and every $N\in\mathbb{N}$, there is a composite polynomial $G_N$ with at most $N$ free parameters such that
\[
\|f-G_N\|_{p}\ \le\ C\,e^{-cN}.
\]
\end{thm}
\footnote{Theorem 2.1 can be modified to work on any finite interval $[A,B]$ using the linear map $L:[A,B]\to [A',B']$ given by $L(x):=\frac{B'-A'}{B-A}x+\frac{AB'-BA'}{A-B}$. We defer this to a future paper.}
\begin{proof}
It suffices to handle a single cusp term $h(|x-a|^\alpha)$ together with the analytic background $H$; we then sum the resulting approximants and take the largest constants. Fix $\alpha=r/s\in(0,1)$ and write $u=t^{1/s}$ for the exact root of $g^s-t=0$.

For notational brevity suppress the $t$--argument. Define
\[
f(g):=g^s-t,\qquad R(g):=f'(g)=s\,g^{s-1},\qquad u:=t^{1/s}\in[0,1],\qquad A:=R(u)=s\,u^{s-1}.
\]
Errors:
\[
e_k:=g_k-u,\qquad \Delta_k:=1-R(g_k)\,y_k.
\]
\emph{Polynomial update identity.} From \eqref{eq:divfree},
\[
R(g_k)\,y_{k+1}=R(g_k)\,y_k\bigl(2-R(g_k)\,y_k\bigr)=2z_k-z_k^2,\qquad z_k:=R(g_k)\,y_k.
\]
Hence the exact identity
\begin{equation}\label{eq:DeltaSquare}
1-R(g_k)\,y_{k+1}=\bigl(1-R(g_k)\,y_k\bigr)^2
\qquad\Longleftrightarrow\qquad
\Delta_{k+1}=\Delta_k^{\,2}.
\end{equation}
\emph{Monotonicity bounds.} Since $g_0=1\ge u$ and $y_0=1/s\in(0,2/R(g_0))$, an induction using \eqref{eq:divfree} and $f(g_k)\ge 0$ shows
\begin{equation}\label{eq:mono}
u\ \le\ g_{k+1}\ \le\ g_k\ \le\ 1,\qquad
0\ <\ R(g_k)\,y_k\ <\ 2,\qquad\text{for all $k\ge 0$.}
\end{equation}
Indeed, $g_{k+1}=g_k-y_{k+1}f(g_k)\le g_k$ and $g_{k+1}\ge u$ because $f(g_k)\ge 0$, $y_{k+1}\ge 0$. Also $0<z_k<2$ holds since $y_{k+1}=y_k(2-R(g_k)y_k)$ preserves the interval $(0,2/R(g_k))$.

\paragraph{Pointwise error recursion for $g_k$.}
By the mean value theorem,
\begin{equation}\label{eq:mvt}
f(g_k)-f(u)=R(\xi_k)\,(g_k-u)=R(\xi_k)\,e_k
\end{equation}
for some $\xi_k$ between $g_k$ and $u$. Using \eqref{eq:divfree},
\[
e_{k+1}=g_{k+1}-u
=g_k-u-y_{k+1}\bigl(f(g_k)-f(u)\bigr)
=\bigl(1-y_{k+1}R(\xi_k)\bigr)\,e_k.
\]
Add and subtract $R(g_k)$:
\[
1-y_{k+1}R(\xi_k)
=\underbrace{\bigl(1-y_{k+1}R(g_k)\bigr)}_{=\Delta_k^{\,2}\ \text{by }\eqref{eq:DeltaSquare}}
\ +\ y_{k+1}\bigl(R(g_k)-R(\xi_k)\bigr).
\]
Hence
\begin{equation}\label{eq:eNextCore}
|e_{k+1}|
\ \le\ \Bigl(\Delta_k^{\,2}+y_{k+1}\,|R(g_k)-R(\xi_k)|\Bigr)\,|e_k|.
\end{equation}
Using $R'(g)=s(s-1)g^{s-2}$ and $g,u\in[0,1]$ we have $|R(g)-R(\xi)|\le s(s-1)\,|g-\xi|\le s(s-1)\,|g-u|=s(s-1)\,|e_k|$. Also $y_{k+1}\le 2/R(g_k)$ by \eqref{eq:mono}. Thus
\begin{equation}\label{eq:eRecTau}
|e_{k+1}|
\ \le\ \Delta_k^{\,2}\,|e_k|\ +\ \frac{2\,s(s-1)}{R(g_k)}\,|e_k|^2.
\end{equation}

\paragraph{Uniform control away from $t=0$.}
Fix any $\tau\in(0,1]$ and restrict to $t\in[\tau,1]$. Then $u=t^{1/s}\ge \tau^{1/s}$ and
\[
R(g_k)\ \ge\ R(u)=A=s\,u^{s-1}\ge s\,\tau^{(s-1)/s}=:A_\tau.
\]
In particular $1/R(g_k)\le 1/A_\tau$. Using this in \eqref{eq:eRecTau} gives
\begin{equation}\label{eq:eRecUniform}
|e_{k+1}|
\ \le\ \Delta_k^{\,2}\,|e_k|\ +\ C_\tau\,|e_k|^2,\qquad
C_\tau:=\frac{2(s-1)}{\tau^{(s-1)/s}}.
\end{equation}
We also need a recursion for $\Delta_k$ that exposes its dependence on $e_{k-1}$. Starting from
\[
\Delta_k
=1-R(g_k)\,y_k
=1-\frac{R(g_k)}{R(g_{k-1})}\,R(g_{k-1})y_{k-1}\bigl(2-R(g_{k-1})y_{k-1}\bigr),
\]
set $\lambda_k:=R(g_k)/R(g_{k-1})\in(0,1]$ (monotonicity \eqref{eq:mono}) and $z_{k-1}:=R(g_{k-1})y_{k-1}$. Then
\[
\Delta_k
=1-\lambda_k\bigl(2z_{k-1}-z_{k-1}^2\bigr)
=1-\lambda_k\bigl(1-\Delta_{k-1}^{\,2}\bigr)
=(1-\lambda_k)+\lambda_k\,\Delta_{k-1}^{\,2}.
\]
Using $1-\lambda_k=(R(g_{k-1})-R(g_k))/R(g_{k-1})$ and the Lipschitz bound on $R$ yields
\begin{equation}\label{eq:DeltaRec}
\Delta_k
\ \le\ \frac{s(s-1)}{R(g_{k-1})}\,|e_{k-1}|\ +\ \Delta_{k-1}^{\,2}
\ \le\ C'_\tau\,|e_{k-1}|\ +\ \Delta_{k-1}^{\,2},\qquad
C'_\tau:=\frac{s-1}{\tau^{(s-1)/s}}.
\end{equation}

\paragraph{Entrance into a quadratic basin with uniform constants.}
Define the combined error $Z_k:=|e_k|+\Delta_k$. From \eqref{eq:eRecUniform} and \eqref{eq:DeltaRec},
\[
\begin{aligned}
|e_{k+1}| &\le \Delta_k^{\,2}|e_k|+C_\tau |e_k|^2
\ \le\ Z_k^3 + C_\tau Z_k^2,\\
\Delta_{k+1} &\le C'_\tau |e_k|+\Delta_k^{\,2}
\ \le\ C'_\tau Z_k + Z_k^2,
\end{aligned}
\]
where we used $\Delta_k\le Z_k$ and $|e_k|\le Z_k$.
Hence there exists $\eta_\tau\in(0,1)$ and a constant $K_\tau>0$ (depending only on $s$ and $\tau$) such that if $Z_k\le \eta_\tau$, then
\[
|e_{k+1}|\le (\eta_\tau + C_\tau) Z_k^2,
\qquad
\Delta_{k+1}\le (C'_\tau \eta_\tau + 1) Z_k^2,
\]
and therefore
\begin{equation}\label{eq:quadRegime}
Z_{k+1} = |e_{k+1}|+\Delta_{k+1}
 \le \bigl(C_\tau + 1 + \eta_\tau(1+C'_\tau)\bigr)\,Z_k^2
  =: K_\tau\,Z_k^2.
\end{equation}
It remains to see that for $t\in[\tau,1]$ the iteration \eqref{eq:divfree} reaches $Z_k\le \eta_\tau$ in a number of steps bounded uniformly in $t$. This follows from the monotonicity \eqref{eq:mono} and continuity: the map $(g,y)\mapsto\bigl(y(2-R(g)y),\ g-y(2-R(g)y)\,f(g)\bigr)$ is continuous on the compact rectangle
\[
\mathcal{R}_\tau:=\bigl\{(g,y):\ u\le g\le 1,\ 0<y\le 2/R(g)\bigr\}\subset (0,1]\times(0,\infty),
\]
and $(g_0,y_0)=(1,1/s)\in\mathcal{R}_\tau$. Since $u\mapsto (u,1/R(u))$ is the unique fixed point of this map on each fiber $t\in[\tau,1]$, standard compactness and the monotone squeezing in \eqref{eq:mono} imply that there exists $k_\tau$ such that $Z_{k_\tau}\le\eta_\tau$ for all $t\in[\tau,1]$. Combining with \eqref{eq:quadRegime} yields
\[
\sup_{t\in[\tau,1]} Z_{k_\tau+\ell}\ \le\ C_1(\tau)\,\rho_\tau^{-\ell}\qquad\text{for some }\rho_\tau>1.
\]
Since $\phi_k(t)=g_k(t)^r$, the map $u\mapsto u^r$ is $r$--Lipschitz on $[0,1]$, and $|u^r-v^r|\le r\,|u-v|$, we obtain
\begin{equation}\label{eq:phiRate}
\sup_{t\in[\tau,1]}\bigl|\phi_k(t)-t^{\alpha}\bigr|
\ \le\ C_2(\tau)\,\rho_\tau^{-k}.
\end{equation}

\paragraph{$L^p$ control near $t=0$.}
On the boundary layer $\{t\le \tau\}$ we simply use $|\phi_k(t)-t^\alpha|\le 1$ and the change of variables $t=|x-a|$ 
to get:
\[
\int_{-1}^1\!\bigl|\phi_k(|x-a|)-|x-a|^\alpha\bigr|^p dx
\ \le\ 2\tau\cdot 1^p\ +\ \int_{\{|x-a|\ge \tau\}}\!\!\!\!\!\!C_2(\tau)^p\,\rho_\tau^{-pk}\,dx
\ \le\ 2\tau + 2\,C_2(\tau)^p\,\rho_\tau^{-pk}.
\]
Taking $(1/p)$th roots gives the $L^p$ bound:
\begin{equation}\label{eq:L2inner}
\bigl\||x-a|^\alpha-\phi_k(|x-a|)\bigr\|_{L^p([-1,1])}
\ \le\ \left({2\tau}\ +\ 2C_2(\tau)^p\,\rho_\tau^{-pk}\right)^{1/p}
\end{equation}
Choosing $\tau=\rho_\tau^{-pk}$ balances the two terms and yields
\begin{equation}\label{eq:innerFinal}
\bigl\||x-a|^\alpha-\phi_k(|x-a|)\bigr\|_{L^p([-1,1])}
\ \le\ C_3\,\rho^{-k}
\quad\text{for some }\rho>1\ \text{independent of }k.
\end{equation}
Here $C_3$ depends only on $p$.
\medskip


\paragraph{Outer analytic approximation and composition.}
Because $h$ is analytic in a complex neighborhood of $[0,1]$, there exist constants $M>0$ and $R>1$ and polynomials $P_m$ of degree $m$ such that
\begin{equation}\label{eq:outerGeom}
\sup_{u\in[0,1]}|h(u)-P_m(u)|\ \le\ M\,R^{-m}.
\end{equation}
This can be proved by expanding $h$ in a Chebyshev series on $[0,1]$ and using Cauchy’s integral estimate for the coefficients on a Bernstein ellipse strictly containing $[0,1]$. For any degree $m$ polynomial $Q_m$ on $[0,1]$, the Markov inequality gives
\begin{equation}\label{eq:markov}
\sup_{u\in[0,1]}|Q_m'(u)|\ \le\ m^2\,\sup_{u\in[0,1]}|Q_m(u)|.
\end{equation}
Now decompose the outer error at $u=t^\alpha$ and $\hat u=\phi_k(t)$:
\[
\bigl|h(u)-P_m(\hat u)\bigr|
\ \le\ \underbrace{|h(u)-P_m(u)|}_{\le M R^{-m}}
\ +\ \underbrace{|P_m(u)-P_m(\hat u)|}_{\le \|P_m'\|_\infty\,|u-\hat u|}.
\]
By \eqref{eq:markov} and \eqref{eq:outerGeom}, $\|P_m'\|_\infty\le m^2\,\|P_m\|_\infty\le m^2\bigl(\|h\|_\infty+M R^{-m}\bigr)\le C_4\,m^2$. Therefore 
\[
\bigl|h(t^\alpha)-P_m(\phi_k(t))\bigr|
\ \le\ M R^{-m}+ C_4\,m^2\,\bigl|\phi_k(t)-t^\alpha\bigr|.
\]
Pull back to $x$ and take $L^p$ quasi-norms.
For $0<p\le 1$ we use the elementary inequality
\begin{equation}\label{eq:subadd_p}
(A+B)^p \le A^p + B^p, \qquad A,B\ge 0,
\end{equation}
which implies the quasi-subadditivity
\(
\|F+G\|_{L^p}^p \le \|F\|_{L^p}^p + \|G\|_{L^p}^p
\)
for measurable $F,G$.
From
\[
|h(|x-a|^\alpha)-P_m(\phi_k(|x-a|))|
\le MR^{-m} + C_4 m^2\,|\phi_k(|x-a|)-|x-a|^\alpha|,
\]
we obtain by \eqref{eq:subadd_p} that
\[
\|h(|x-a|^\alpha)-P_m(\phi_k(|x-a|))\|_{L^p}^p
\le (MR^{-m})^p\,\|1\|_{L^p}^p
+ (C_4 m^2)^p\,\|\phi_k(|x-a|)-|x-a|^\alpha\|_{L^p}^p.
\]
Since $\|1\|_{L^p([-1,1])}^p=\int_{-1}^1 1\,dx=2$ and using \eqref{eq:innerFinal}, we conclude
\begin{equation}\label{eq:15p}
\|h(|x-a|^\alpha)-P_m(\phi_k(|x-a|))\|_{L^p}
\le \Bigl( 2(MR^{-m})^p + (C_5 m^2 \rho^{-k})^p \Bigr)^{1/p}.
\end{equation}
Moreover, for any $p>0$ and $a,b\ge0$ we have $(a+b)^{1/p}\le a^{1/p}+b^{1/p}$, hence
\[
\Bigl( 2(MR^{-m})^p + (C_5 m^2 \rho^{-k})^p \Bigr)^{1/p}
\le 2^{1/p}MR^{-m}+C_5m^2\rho^{-k}.
\]
Defining
\[
M' := 2^{1/p}M,
\]
which may be large when $p$ is small but is a fixed constant for fixed $p$, yields
\begin{equation}\label{eq:outerInnerCombined}
\|h(|x-a|^\alpha)-P_m(\phi_k(|x-a|))\|_{L^p([-1,1])}
\le M'R^{-m} + C_5 m^2\rho^{-k}.
\end{equation}
For $1\le p<\infty$, the argument simplifies, since Minkowski’s inequality holds in $L^p$ and
yields the same conclusions with constants depending only on $p$.

\paragraph{Adding the analytic background and balancing parameters.}
Approximate $H$ by a degree $m$ polynomial $H_m$ with the same geometric rate as in \eqref{eq:outerGeom}. Then
\[
\|H-H_m\|_{L^p([-1,1])}\ \le\ C_6\,R^{-m}.
\]
Here, $C_6$ depends only on $p$.
\medskip

Define the deep approximant
\[
G_{m,k}(x):=H_m(x)+P_m\!\left(\phi_k(|x-a|)\right).
\]
Combining the bounds for $H-H_m$ and \eqref{eq:15p} using \eqref{eq:subadd_p} (for $0<p\le1$) gives
\[
\|H+h(|x-a|^\alpha)-G_{m,k}\|_{L^p}^p
\le \|H-H_m\|_{L^p}^p + \|h(|x-a|^\alpha)-P_m(\phi_k(|x-a|))\|_{L^p}^p.
\]
Using $\|H-H_m\|_{L^p}\le C_6R^{-m}$ and \eqref{eq:15p}, we obtain
\[
\|H+h(|x-a|^\alpha)-G_{m,k}\|_{L^p}
\le \Bigl( (C_6R^{-m})^p + 2(MR^{-m})^p + (C_5 m^2\rho^{-k})^p \Bigr)^{1/p}.
\]
Absorbing constants yields
\[
\|H+h(|x-a|^\alpha)-G_{m,k}\|_{L^p([-1,1])}
\le C_7R^{-m} + C_5 m^2\rho^{-k}.
\]
Here, $C_7$ depends only on $p$.
For a desired parameter budget $N$ choose $m=\lfloor \gamma k\rfloor$ with a fixed $\gamma>0$ so that $N\asymp m+k$. Then
\[
\|f-G_{m,k}\|_{L^p([-1,1])}
\ \le\ C_8\Bigl(R^{-\gamma k}+k^2\,\rho^{-k}\Bigr)
\ \le\ C_9\,e^{-cN}
\]
for suitable $c\in(0,\min\{\log R^\gamma,\log \rho\}/(1+\gamma))$, absorbing polynomial factors of $k$ into the constant. Finally, summing over finitely many cusp terms $j=1,\dots,M$ only changes the leading constant. Here, both 
$C_8,C_9$ depend only on $p$. This proves the theorem.
\end{proof}

\paragraph{Discussion of rates and parameters.}
The inner iteration \eqref{eq:divfree} produces a chain of bounded--degree “layers,” so its coefficient count grows linearly with $k$. The outer fit adds $m+1$ coefficients. Thus $N\asymp m+k$, and the two–term error in \eqref{eq:outerInnerCombined} is balanced by taking $m\propto k$, which yields an overall $L^2$ error that decays exponentially in $N$. The key algebraic mechanisms are the exact identity \eqref{eq:DeltaSquare} for the reciprocal update and the mean–value reduction \eqref{eq:mvt} for the primal update, which together force a quadratic basin with uniform constants away from $t=0$ and a diminishing boundary layer near $t=0$ whose measure can be matched to the inner rate.

\begin{figure}[h]
    \centering
    \begin{minipage}{0.50\textwidth} 
        \centering
        \includegraphics[width=\linewidth]{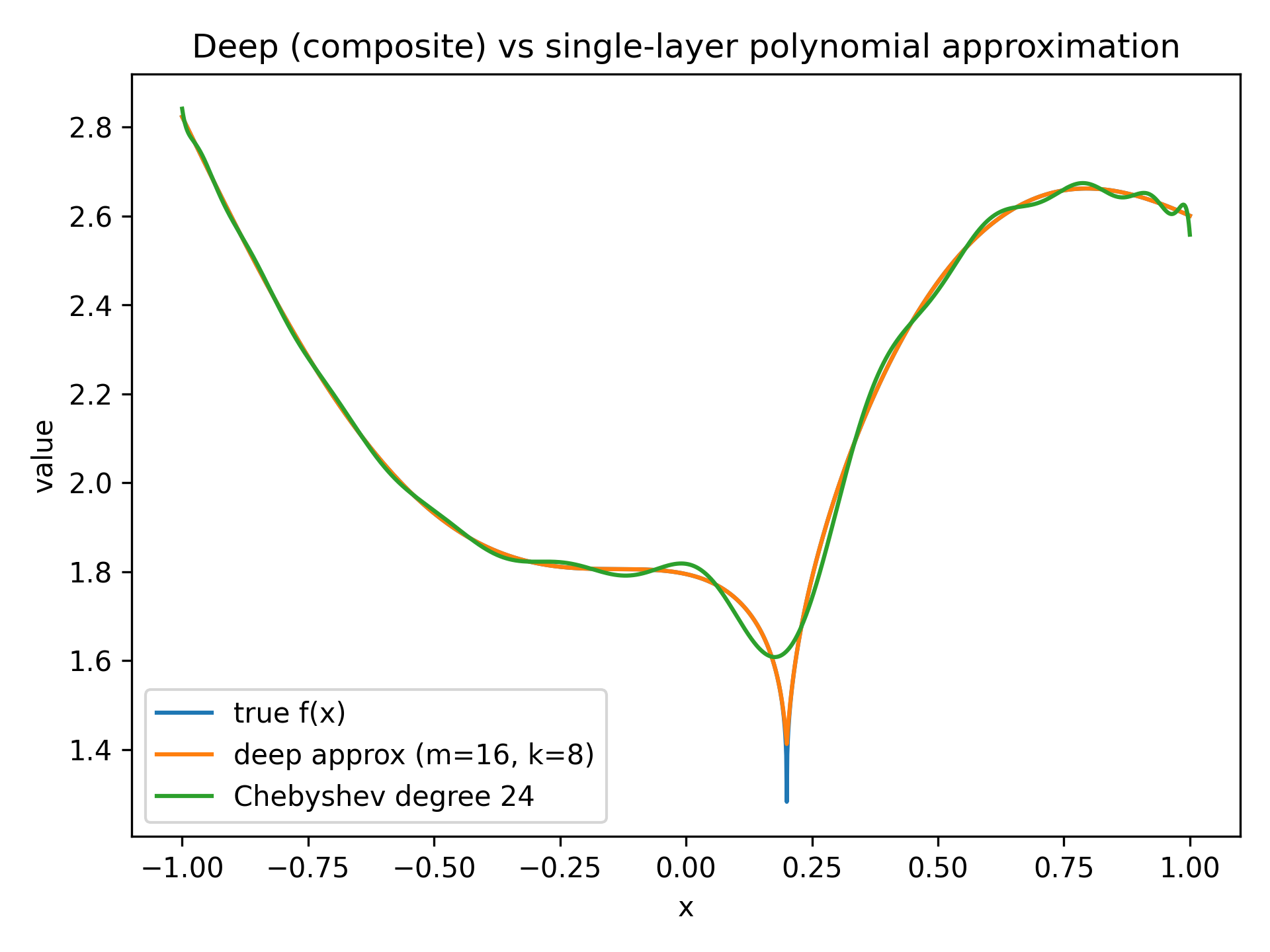}
        \caption*{(a) Approximation of cusp type function}
    \end{minipage}\hfill
    \begin{minipage}{0.50\textwidth}
        \centering
        \includegraphics[width=\linewidth]{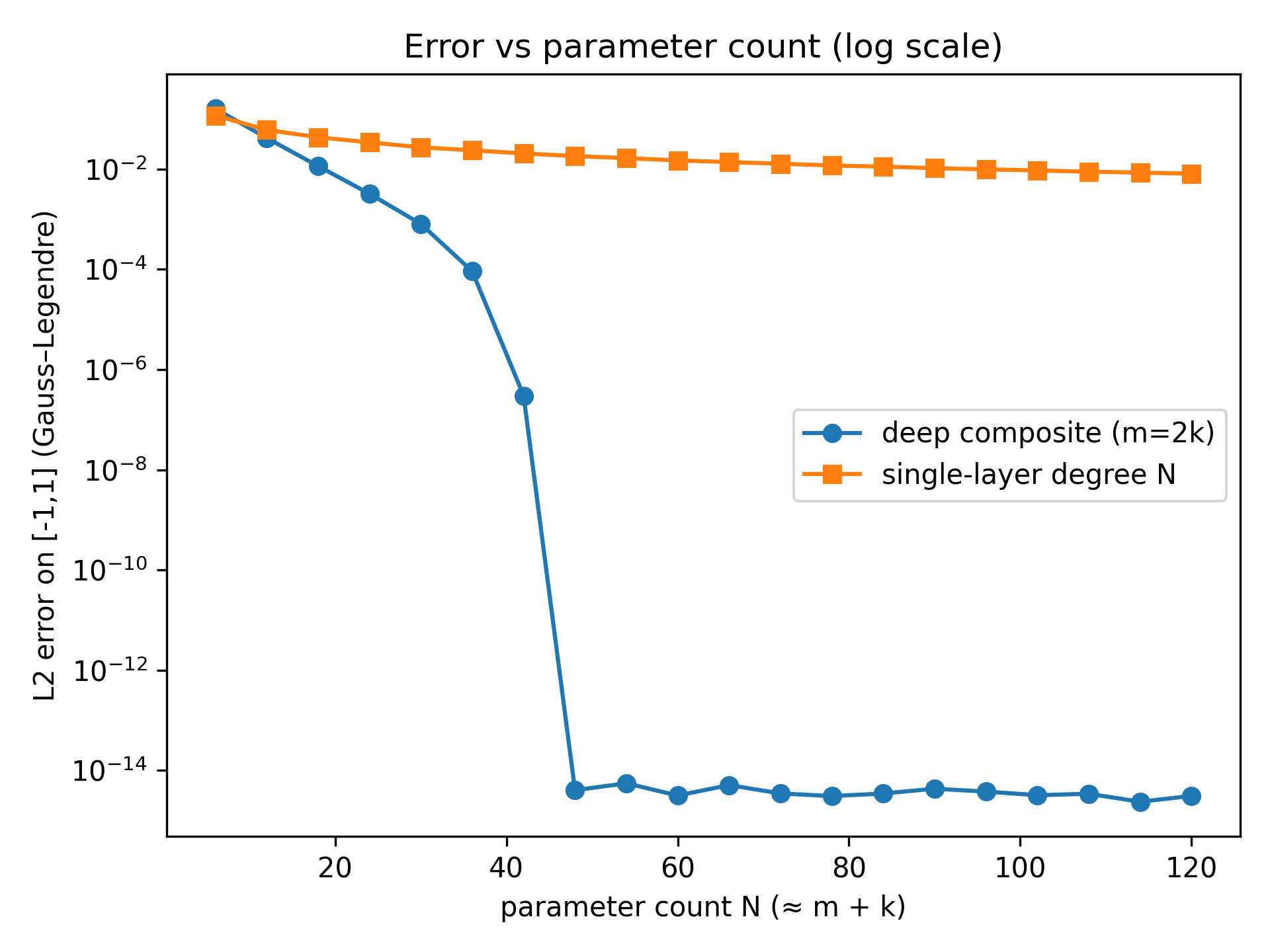}
        \caption*{(b) Pointwise absolute Error (log scale)}
    \end{minipage}
    \caption{Comparison of approximants for cusp type function}
    \label{fig:curspapprox}
\end{figure}

\begin{figure}[h]
    \centering
    \begin{minipage}{0.50\textwidth}
        \centering
        \includegraphics[width=\linewidth]{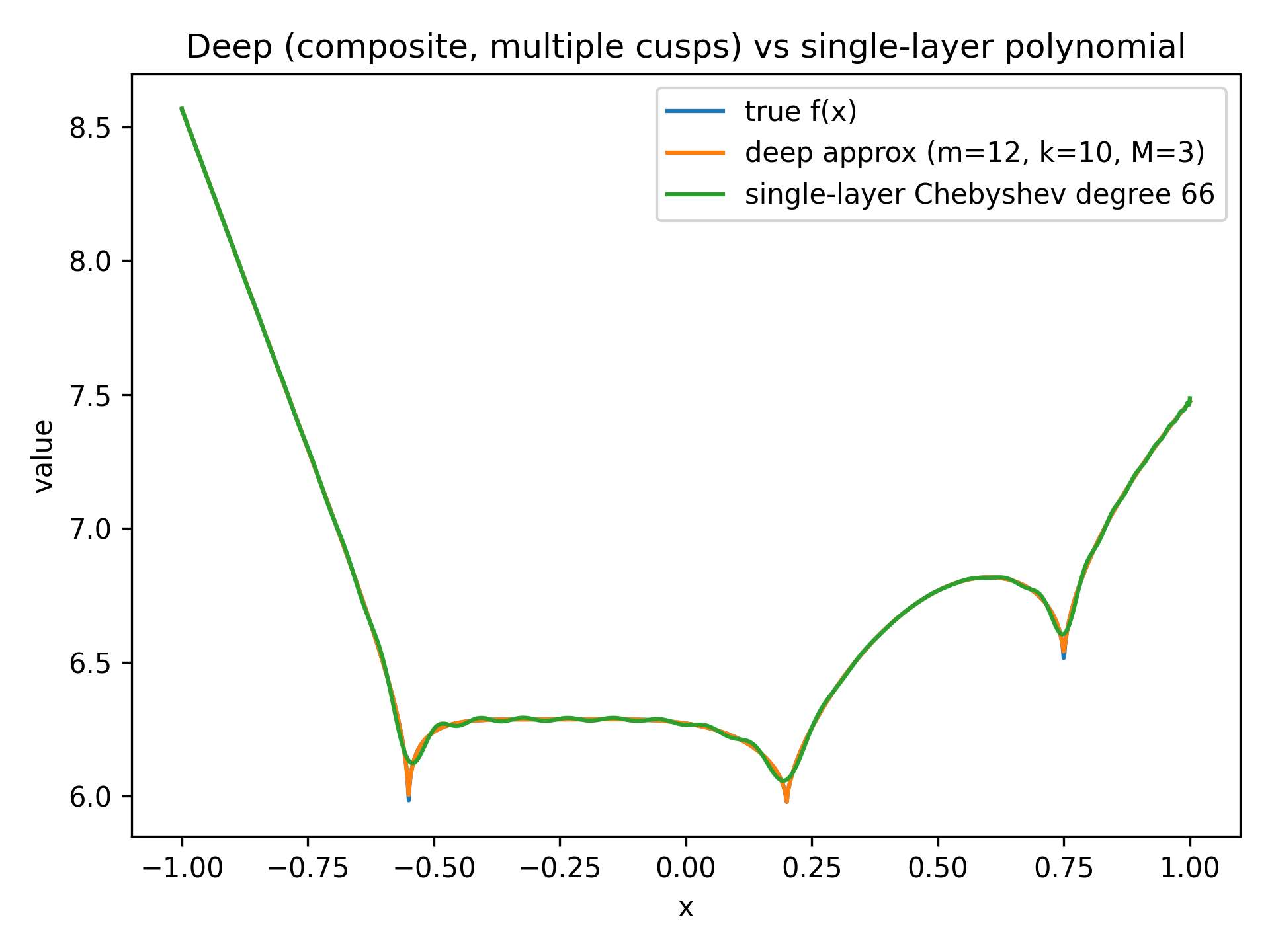}
        \caption*{(a) Approximation of multiple cusp type function}
    \end{minipage}\hfill
    \begin{minipage}{0.50\textwidth}
        \centering
        \includegraphics[width=\linewidth]{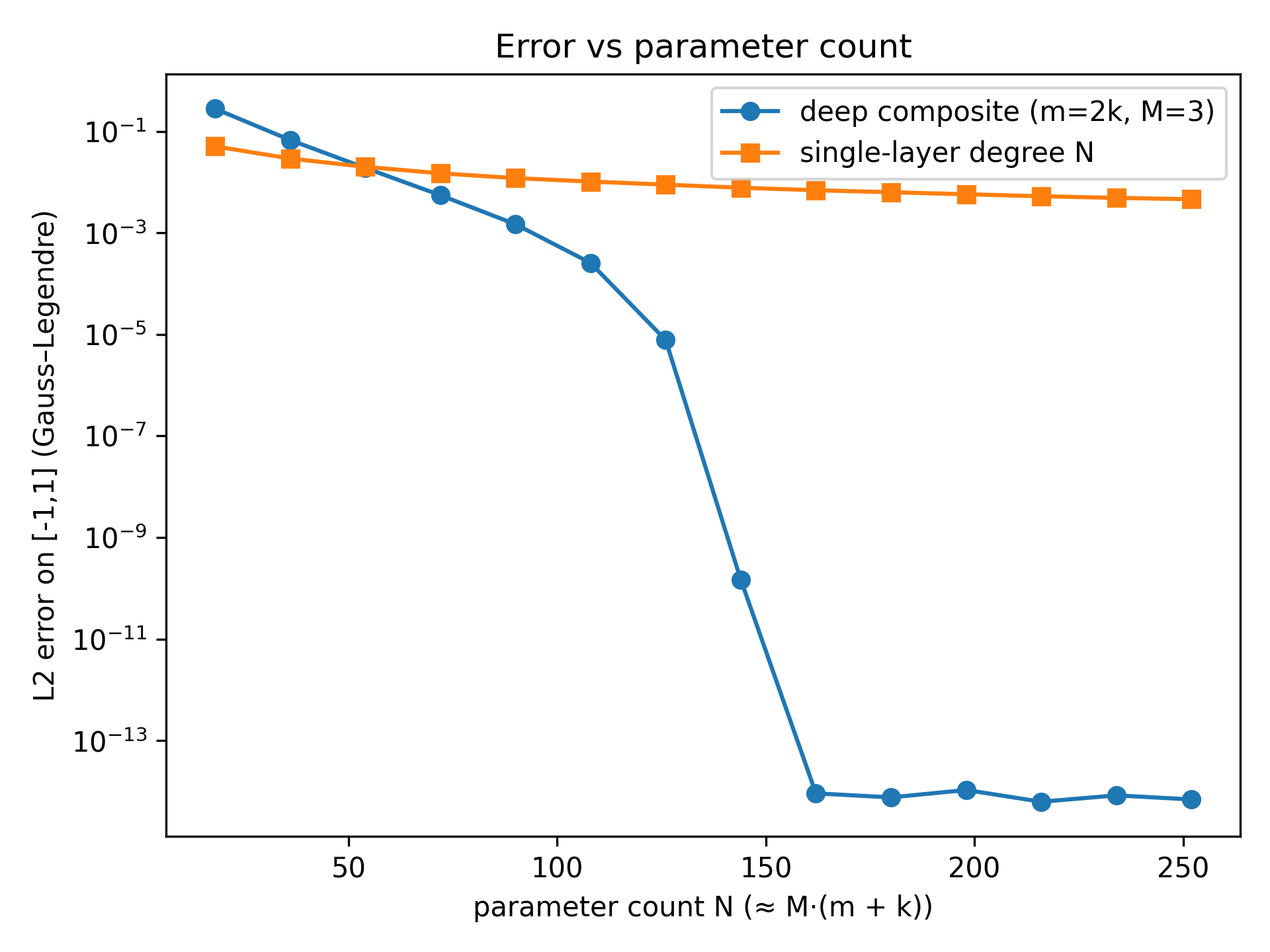}
        \caption*{(b) Pointwise absolute Error (log scale)}
    \end{minipage}
    \caption{Comparison of approximants for multiple cusp type function}
    \label{fig:mutiplecurspapprox}
\end{figure}

Figures~\ref{fig:curspapprox} and~\ref{fig:mutiplecurspapprox} illustrate the approximation behavior of the proposed deep composite polynomial model for functions with algebraic cusp singularities.
In Figure~\ref{fig:curspapprox}, we consider a single cusp located at $a=0.2$ with exponent $\alpha=r/s=1/3$, where the target function has the form
$f(x)=H(x)+h(|x-a|^\alpha)$ with an analytic envelope $h$.
The left panel shows a representative approximation using parameters $(m,k)$, while the right panel plots the $L^2([-1,1])$ error, computed via Gauss--Legendre quadrature, as a function of the effective parameter count $N\approx m+k$.
The deep composite approximation exhibits a markedly faster error decay compared to a single-layer Chebyshev polynomial of comparable degree.
Figure~\ref{fig:mutiplecurspapprox} shows analogous results for a function with multiple cusps at distinct locations and exponents, where the total parameter count scales as $N\approx M(m+k)$ with $M$ cusp terms.
In both cases, the empirical results confirm the exponential-in-$N$ convergence predicted by the theory and demonstrate a clear advantage of composite polynomial constructions over single-layer polynomial approximations.

We illustrate the approximation properties of the composite polynomial model on two–dimensional test functions with angular cusp singularities.
Let $(x,y)\in[-1,1]^2$, define polar coordinates
\[
r = \sqrt{x^2+y^2}, \qquad 
\theta = \operatorname{atan2}(y,x)\in(-\pi,\pi],
\]
and construct a star–shaped domain through a radial profile $R_\star(\theta)$.
Specifically, we consider
\[
R_{\star}(\theta)
=
R_{0}
+
\sum_{j=0}^{K_{\star}-1}
W_{j}
\exp\!\bigl(-\lambda_{j}\,\lvert \theta - \theta_{j} \rvert^{\alpha_j}\bigr),
\]
where $K_\star$ denotes the number of star tips, $\theta_j$ are their angular locations, $W_j>0$ their amplitudes, $\lambda_j>0$ decay rates, and $\alpha_j\in(0,1)$ control the sharpness of the angular cusps.
The associated level–set function is
\[
f(x,y)=\tanh\!\bigl(\gamma\,(R_{\star}(\theta)-r)\bigr),
\]
so that the zero contour $f=0$ traces a simply connected star whose boundary is smooth away from the tips and exhibits algebraic angular cusps at $\theta=\theta_j$.

\begin{figure}[h]
    \centering
    \includegraphics[width=1.0\textwidth]{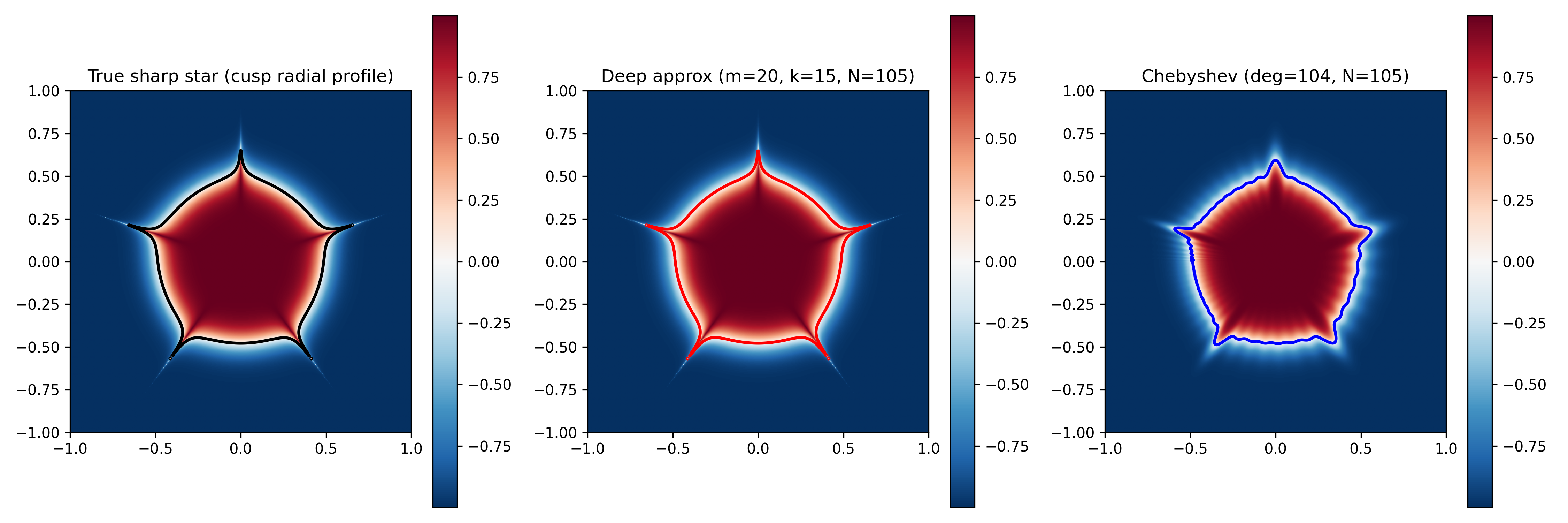}
    \caption{
    Sharp five–point star with algebraic angular cusps.
    Left: true level set $f(x,y)=\tanh(\gamma(R_\star(\theta)-r))$ (black); middle: deep composite approximation ($K_\star=5$, $m=20$, $k=15$, $N=105$, red); right: Chebyshev baseline of degree $D=104$ with matched parameter count (blue), with $L^2$ errors $3.17\times10^{-3}$ (deep) versus $3.55\times10^{-2}$ on a $400\times400$ grid.
    }
    \label{fig:cusp2d}
\end{figure}

We first consider a symmetric five–point star with $K_\star=5$, equally spaced spike locations $\theta_j=\theta_0+2\pi j/5$, and a fixed cusp exponent $\alpha_j=1/3$ for all $j$ (Figure~\ref{fig:cusp2d}, left).
We approximate $R_\star(\theta)$ using the proposed deep composite architecture, which applies an inner division–free Newton map $\phi_k$ to approximate the fractional powers $\lvert\theta-\theta_j\rvert^{\alpha_j}$ and then fits outer Chebyshev polynomials of degree $m$ to the resulting envelopes.
As a baseline, we fit a single one–dimensional Chebyshev polynomial in $\theta$ with the same total number of scalar parameters and lift it back to two dimensions using the same level–set construction.
With $m=20$ and $k=15$, both methods use $N=K_\star(m+1)=105$ parameters.
On a $400\times400$ grid over $[-1,1]^2$, the deep composite approximation attains an $L^2$ error of $3.17\times10^{-3}$, compared to $3.55\times10^{-2}$ for the Chebyshev baseline, while more accurately resolving the sharp cusp geometry (Figure~\ref{fig:cusp2d}).

We next consider an uneven eight–point star with $K_\star=8$, where the spike locations $\theta_j$ are randomly perturbed from uniform spacing and the parameters $\alpha_j$, $W_j$, and $\lambda_j$ vary across tips, producing nonuniform cusp sharpness and asymmetric geometry (Figure~\ref{fig:cusp2d_uneven8}).
Using $m=22$ and $k=16$, corresponding to $N=184$ parameters, the deep composite model achieves an $L^2$ error of $9.74\times10^{-3}$ on a $420\times420$ grid, compared to $2.40\times10^{-2}$ for the Chebyshev approximation.
Despite the increased geometric irregularity, the composite model consistently captures the localized cusp behavior more accurately than the global polynomial baseline.

\begin{figure}[h]
    \centering
    \includegraphics[width=1.0\textwidth]{cusp_2d_uneven_8point.png}
    \caption{
        Uneven eight–point star with algebraic angular cusps.
        Left: true level set; middle: deep composite approximation ($K_\star=8$, $m=22$, $k=16$, $N=184$); right: Chebyshev baseline with matched parameter count, with zero contours shown in black, red, and blue.
    }
    \label{fig:cusp2d_uneven8}
\end{figure}


\bibliographystyle{elsarticle-num} 
\bibliography{refs}







\end{document}